# A Nationwide Multi-Location Multi-Resource Stochastic Programming Based Energy Planning Framework

Razan A. H. Al-Lawati[1], Tasnim Ibn Faiz[1], Md. Noor-E-Alam[1*]

[1]Dept. of Mechanical and Industrial Engineering, Northeastern University, 360 Huntington Avenue, 02115, Boston, MA, USA

* Corresponding author Tel.: +1-617-373-2275

E-mail address: mnalam@neu.edu

**Abstract**

With the escalating global energy demands, nations face greater urgency in diversifying their electricity market portfolio with various energy sources. The efficient utilization of these fuel resources necessitates the integrated optimization of the entire energy portfolio. However, the inherent uncertainty of variable renewable energy generators renders traditional deterministic planning models ineffective. Conversely, expansive stochastic models that capture all the complexities of the power generation problem are often computationally intractable. This paper introduces a scalable energy planning framework tailored for policymakers, which is adaptable to diverse geopolitical contexts, including individual countries, regions, or coalitions formed through energy trade agreements. A two-stage stochastic programming approach, including a scenario-based Benders decomposition method, is employed to achieve computational efficiency. Scenarios for representing parameter uncertainty are developed using a *k*-means clustering algorithm. An illustrative case study is presented to show the practical application of the proposed framework and its ability to inform policymaking by identifying actionable cost-reduction strategies. The findings underscore the importance of promoting resource coordination and demonstrate the impact of increasing interchange and renewable energy capacity on overall gains. Notably, the study shows the superiority of stochastic optimization over deterministic models in energy planning.

## 1. Motivation for an integrated model

The global energy demand continues to increase [1] for all fuel sources [2], while the supply of fossil fuels continues to deplete and cause adverse environmental impacts [3]. As such, evaluating how best to utilize alternative fuel sources is imperative. Renewable energy development has become a joint effort of all countries worldwide [4]. In 2022, renewable energy resources comprised 40% of the world's installed power capacity, which is likely to grow [5]. However, renewable energy production can be highly variable, and the uncertainty associated with their energy production makes forecasting and decision-making difficult for energy regulators [6].

### 1.1. Optimizing regional energy resources

It is essential to optimize all energy resources concurrently within a single framework to improve the efficient utilization of the available energy resources. However, a gap exists in the literature aimed at achieving optimal designs of hybrid energy systems, which was highlighted by the excellent review work of Zahraee et al. [7]. Out of the over one hundred research works they reviewed, they only found three models that incorporated multiple sources of conventional and renewable energy within one framework. The first two models addressed Generation Expansion Planning (GEP) problems and employed multi-objective optimization approaches [8], [9], whereas the third utilized an algorithm that combines particle swarm optimization and a genetic algorithm to forecast energy demand [10].

Subsequent works focused on hybrid energy systems, including Ghaithan et al. [11], which proposes a grid-connected solar-wind system to power a small-scale Reverse Osmosis (RO) desalination unit in Saudi Arabia. Abdelshafy et al. [12] optimized a grid-connected double storage system in Egypt consisting of pumped-storage hydropower (PSH) and batteries supplied by photovoltaics (PV) and wind turbines. Micangeli et al. [13] optimized a mini-grid for remote areas in Kenya with an off-grid hybrid power station based on wind, PV-solar, and diesel generation. Alayi et al. [14] optimized electricity supply to remote areas in Iran using hybrid energy systems. Mallek et al. [15] used HOMER simulation software to optimize several hybrid energy system models consisting of solar PV, wind, and the national grid in a desalination plant in Tunisia. Kefale et al. [16] optimized the sizing and siting of the solar photovoltaic system in the Ethiopian radial distribution system using selective particle swarm optimization (SPSO). Lastly, Hassan et al. [17] developed a simulation to optimize two PV power system configurations to encourage residential PV power systems for household usage and feed the national grid of Iraq.

While each of these models presents valuable insights into optimizing energy systems, they typically do not account for the full spectrum of energy sources concurrently or at a scale that encompasses national energy decision-making. Our model advances the field by integrating an array of both conventional and renewable energy sources within a single framework tailored for national-level optimization and policy development.

### 1.2. Optimizing the sharing of resources

Globally, some economies are becoming more energy-independent, while others depend more heavily on imports to satisfy their energy needs [18]. Cooperation and sharing of energy resources are catalysts for sustainable growth. It strengthens national energy security, reduces the cost of supplies, and minimizes the impact of price volatility [19]. International energy trade has been studied using networks to understand their relationships and extract patterns, such as the work of Xu et al. [20] that quantifies the characteristics of the complex relationship of international energy trade through network reconstruction approaches, or Berdysheva et al. [18] that examines the evolution of the global network. There has also been work looking at embodied energy, which is the total direct and indirect energy required to produce goods and services incorporated in the product, such as the work of Chen et al. [21]. Other work looks at the impact of different factors on energy trade. Li et al. [22] developed a regression-discontinuity model to study the effect of geopolitics.

The existing literature on international energy trade and cooperation limits its focus on studying global characteristics and trends rather than optimizing direct energy trade between regions. There remains a significant need for optimization models that can operationalize energy trade and cooperation concepts in national energy planning. Our study presents a model and an optimization framework that nations can leverage to make energy resource-sharing decisions and enhance energy security, economic efficiency, and sustainability.

### 1.3. The use of decomposition techniques

Parameter uncertainty presents a non-trivial challenge to decision-making. Stochastic programming is a popular approach for handling uncertainty in developing decision-making frameworks since it allows uncertainty to be incorporated into the model [23]. A common technique for managing parameter uncertainty in stochastic programming is introducing many scenarios, each representing a realization of the uncertain parameters. However, the problem becomes computationally intractable when the number of scenarios is too large. To that end, several decomposition algorithms have been developed to gain computational efficiency.

Benders decomposition, a method that facilitates solving large-scale optimization problems by partitioning the problem into multiple small-size problems [24], is increasingly popular in energy applications. In the realm of unit commitment and scheduling, Benders decomposition has been effectively used. Lopez-Salgado et al. [25] and Saberi et al. [26] applied it to the unit commitment problem, focusing on power generator scheduling. Additionally, Xia et al. [27] used its bi-level nature to preserve the market participants' private information while optimizing the scheduling of energy exchange. Li et al. [28] utilize a solution method based on Benders Decomposition and Modified stochastic Dual Dynamic Programming (MSDDP) [29] to obtain cost-optimal power generation plans for hydropower producers. In the context of smart grids and energy hubs, Soares et al. [30] optimized resource planning and operation from the perspective of aggregators, while Mansouri et al. [31] applied similar techniques for energy hub optimization. Regarding energy storage and investment decisions, Abdulgalil et al. [32] optimized the sizing of energy storage systems (ESS) under wind uncertainties. Baringo et al. [33] tackled an investment portfolio problem, focusing on contract management and site selection for building wind power facilities, respectively. Lastly, Jenabi et al. [34] utilized Benders decomposition to solve an integrated resource planning problem for power systems expansion.

### 1.4. Contribution

Despite energy planning being extensively researched in literature and Benders decomposition and proving its practicality for solving large-scale problems, no models have been developed to optimize the nationwide energy decision problem.

The main contributions of this work are as follows:

1. A decision problem, from the perspective of a policymaker, is considered, which highlights the need for holistic models that capture a nation's diverse energy mix and address the strategic interplay between different fuel sources. By doing so, the model provides a platform for identifying systemic optimization strategies for resource utilization and cost efficiency.
2. The study employs a two-stage stochastic programming approach with scenarios generated via a k-means clustering algorithm to characterize uncertainty. A scenario-based Benders decomposition approach is utilized to gain computational efficiency for a large-scale model.
3. The proposed decision framework is generalizable to any country, region, or any group of countries participating in energy trade agreements to gain a competitive advantage.
4. A case study is presented to illustrate the practical application of the framework. We extract managerial insights and offer evidence-based recommendations through numerical experiments, validating the proposed methods and highlighting their potential for policy intervention and strategic planning in the energy sector.

Our work aims to overcome the shortcomings identified in the current literature by providing a comprehensive and scalable model for national energy planning. Our approach not only enhances the methodological toolkit for policymakers but also sets a new precedent in the energy optimization field for handling the intricacies of modern, diversified energy portfolios.

## 2. Nomenclature

| Sets and subindices | |
|---|---|
| $\iota$ | Iteration |
| $S$ | The set of all scenarios under consideration |
| $s$ | Subindex for scenarios, $s \in S$ |
| $F$ | The set of all fuel sources under consideration |
| $F_C$ | The set of all fuel sources that have a constant generation rate over the entire time period |
| $F_U$ | The set of all fuel sources whose generation rate can be controlled |
| $F_V$ | The set of all VRRGs |
| $f$ | The subindex for fuel sources, $f \in F$ |
| $L$ | The set of all regions/locations |
| $l_i, l_j$ | Subindices for region, $l_i, l_j \in L$ |

| Parameters | |
|---|---|
| $\bar{P}_{l_i,f}$ | Rated power at location $l_i$ for generator type $f$ (MWh) |
| $\hat{P}_{l_i,f}$ | Available power at location $l_i$ for generator type $f$ (MWh) |
| $\hat{P}_{l_i,f,s}$ | Forecasted power available at location $l_i$ for generator type $f$ in scenario $s$ (MWh) |
| $\bar{T}_{l_i,l_j}$ | The capacity of transmission from the location $l_i$ to the location $l_j$ (MWh) |
| $\rho_s$ | Probability of scenario $s$ |
| $D_{l_i,s}$ | Forecasted demand at location $l_i$ in scenario $s$ (MWh) |
| $C^{trans}_{l_i,l_j}$ | Cost of energy interchange from location $l_i$ to location $l_j$ ($/MWh) |
| $C^{prod}_{l_i,f}$ | Cost of generating energy for generator type $f$ at location $l_i$ ($/MWh) |
| $\gamma^{short}_{l_i}$ | Cost of an energy shortage ($/MWh) |
| $\gamma^{trans}_{l_i,l_j}$ | Penalty for deviating from planned interchange amount ($/MWh) |

| Decision variables | |
|---|---|
| $P^{prod}_{l_i,f,s}$ | Amount of energy produced at location $l_i$ for generator type $f$ in scenario $s$ (MWh) |

| $P_{l_i,l_j}^{trans,plan}$ | Amount of energy pre-allocated for interchange from location $l_i$ to location $l_j$ (MWh) |
| --- | --- |
| $P_{l_i,l_j,s}^{trans,actual}$ | Amount of energy used in interchange from location $l_i$ to location $l_j$ in scenario $s$ (MWh) |
| $P_{l_i,l_j,s}^{trans,dev}$ | The deviation between the amount of energy allocated and used in interchange from location $l_i$ to location $l_j$ in scenario $s$ (MWh) |
| $P_{l_i,s}^{short}$ | Energy shortage at location $l_i$ in scenario $s$ (MWh) |
| $P_{l_i,s}^{excess}$ | Energy excess at location $l_i$ in scenario $s$ (MWh) |

Dual variables

| $\mu_s^{trans,plan(\iota)}$ | Sensitivity for first-stage interchange allocation variable |
| --- | --- |

Benders cut variables

| $\alpha^{down}$ | A very large negative constant to avoid an unbounded solution in the first iteration |
| --- | --- |
| $\alpha^{(\iota)}$ | An auxiliary variable representing the objective function of the subproblems |

## 3. Modeling Assumptions

The model is constructed for a national electric grid from the perspective of a decision maker, which reflects the complexities and opportunities for coordination across different fuel types and regional energy systems. The adaptable model can accommodate various configurations, including single or multiple fuel types and regions.

Energy demand loads are assumed to be inelastic to price, which emanates from the fact that most consumers have limited exposure to real-time electricity prices and thus have little incentive to modify their consumption [35], [36]. The model treats generation, transmission, and shortage costs as deterministic parameters.

The unit commitment (UC) problem, which involves scheduling generators to meet demand while adhering to operational constraints and minimizing costs, is acknowledged [37]. However, given the national scope of our study, we simplify by assuming that differences between generators of the same type are marginal; thus, UC constraints are not explicitly accounted for in the model formulation.

Although Energy Storage Systems (ESSs) play a more significant role in the functionality of power grids as the installed capacity continues to grow [38], the current total capacity is still fractional compared to the overall generating capacity. In 2022, the global rated power of the EES was only 174 GW [39], compared to an installed generation capacity of more than 8 TW [40]. Consequently, in this work, ESSs are not included in the decision-making process, allowing decisions for each time step to be considered independently, simplifying the model and enhancing the computational efficiency.

## 4. Model Formulation

The purpose of this model is to assist decision-makers in making long-term strategic decisions for energy systems. The energy landscape is inherently characterized by uncertainty, with fluctuating energy availability, demand projections, and volatile market conditions, and decisions are often made before many parameters are known for certain. The presented model is designed to integrate uncertainty directly into the decision-making process for the planning of regional energy interchange. Using stochastic programming and scenario-based Benders decomposition, the proposed framework addresses the uncertainty and sets a benchmark for modeling regional energy planning at a scale that reflects the complex interdependencies and dynamics of contemporary energy systems.

### 4.1. Two-stage stochastic programming model

Two-stage stochastic programming involves dividing decisions into first and second-stage variables based on their timing relative to the realization of uncertain data. Scenarios are used to represent various realizations of uncertain parameters, and each scenario has its probability of occurrence. First-stage variables are independent of these scenarios, while second-stage variables are defined over each scenario. In this work, the amount of energy allocated for regional interchange, denoted as $P_{l_i,l_j}^{trans,plan}$, is a first-stage variable. Conversely, decisions made during energy market participation are categorized as second-stage variables. Equations (1)-(10) detail the two-stage stochastic programming model.

$$\underset{\substack{P_{l_i,f,s}^{prod}, P_{l_i,l_j}^{trans,plan}, P_{l_i,l_j,s}^{trans,actual} \\ P_{l_i,l_j,s}^{trans,dev}, P_{l_i,s}^{short}, P_{l_i,s}^{excess}}}{\text{Minimize}} \sum_{l_i \in L} \sum_{l_j \in L} \sum_{f \in F} \sum_{s \in S} \left( C_{l_i,l_j}^{trans} P_{l_i,l_j}^{trans,plan} \right.$$

$$\left. + \rho_s \left( C_{l_i,f}^{prod} P_{l_i,f,s}^{prod} + \gamma_{l_i}^{short} P_{l_i,f,s}^{short} + \gamma_{l_i,l_j}^{trans} P_{l_i,l_j,s}^{trans,dev} \right) \right) \quad (1)$$

Subject to

$$P_{l_i,f,s}^{prod} \leq \bar{P}_{l_i,f} \quad \forall l_i \in L, \forall f \in F, \forall s \in S \quad (2)$$

$$P_{l_i,f,s}^{prod} = \hat{P}_{l_i,f} \quad \forall l_i \in L, \forall f \in F_C, \forall s \in S \quad (3)$$

$$P_{l_i,f,s}^{prod} \leq \hat{P}_{l_i,f} \quad \forall l_i \in L, \forall f \in F_U, \forall s \in S \quad (4)$$

$$P_{l_i,f,s}^{prod} \leq \hat{P}_{l_i,f,s} \quad \forall l_i \in L, \forall f \in F_V, \forall s \in S \quad (5)$$

$$P_{l_i,l_j}^{trans,plan} \leq \bar{T}_{l_i,l_j} \forall l_i \in L, \forall l_j \in L \tag{6}$$

$$P_{l_i,l_j,s}^{trans,actual} \leq P_{l_j,l_i}^{trans,plan} \quad \forall l_i \in L, \forall l_j \in L, \forall s \in S \tag{7}$$

$$P_{l_i,l_j,s}^{trans,dev} = P_{l_j,l_i}^{trans,plan} - P_{l_i,l_j,s}^{trans,actual} \quad \forall l_i \in L, \forall l_j \in L, \forall s \in S \tag{8}$$

$$\sum_{f \in F} P_{l,f,s}^{prod} - \sum_{l_j \in L} P_{l_i,l_j,s}^{trans,actual} + \sum_{l_j \in L} P_{l_j,l_i,s}^{trans,actual} + P_{l_i,s}^{short} - P_{l_i,s}^{excess} = D_{l_i,s} \quad \forall s \in S, \forall l_i \in L \tag{9}$$

$$P_{l_i,f,s}^{prod}, P_{l_j,l_i}^{trans,plan}, P_{l_i,l_j,s}^{trans,actual}, P_{l_i,s}^{short}, P_{l_i,s}^{excess} \geq 0 \quad \forall l_i, l_j \in L, \forall f \in F, \forall s \in S \tag{10}$$

The objective function of the stochastic optimization problem (1) is to minimize the expected total cost, including energy production and transmission costs and penalties for not meeting demand or committed transmission amounts. Since the amount of energy produced and the resulting deviations from the demand and the committed transfers are scenario-dependent, the related costs are probability-weighted.

Constraint (2) restricts energy production to the generators' capacity, while constraints (3)-(5) limit it based on available energy. The model categorizes generators into three types: inflexible generators like nuclear power (fixed energy output per time step), resource-dependent generators like petroleum (energy output based on available resources), and VRRGs, whose output depends on unpredictable environmental factors. These generator types are represented as $F_C$ (inflexible), $F_U$ (resource dependent), and $F_V$ (VRRGs), with the collective set $F$ encompassing all fuel types. Constraint (6) governs the energy transfer between locations, capped by the regional interchange capacities. As illustrated in Figure 2, not all locations are interconnected. The model reflects this by setting the energy transfer between disconnected locations to zero in constraint (6), emphasizing realistic network limitations.

Constraint (7) limits energy transfer to the pre-allocated amount for interchange, while constraint (8) calculates the differential between allocated and actual transfer amounts. Constraint (9) acts as a balancing constraint, accounting for energy shortages or surpluses. Lastly, constraint (10) imposes non-negativity restrictions, ensuring all calculated values are positive or zero.

### 4.2. Decomposed model

Benders decomposition is employed to improve the computational efficiency of the solution process of the resulting large-scale problem with numerous scenarios. With Benders decomposition,

'complicating variables' or 'coupling variables' refer to those decision variables that divide the problem into multiple, independently solvable subproblems if set to a constant. A master problem is constructed to ascertain suitable values for these complicating variables. As illustrated in Figure 1, the algorithm runs iteratively.

Once the master problem sets the values for the complicating variables, these fixed values are used to solve the subproblems. After solving the subproblems, Benders cuts are generated using the dual information from the constraints involving the complicating variables. The master problem is then re-solved with these Benders cuts, narrowing the feasible region, and this process repeats until it converges to an optimal solution. Benders decomposition guarantees optimality when the original problem's objective function is convex [41], [42].

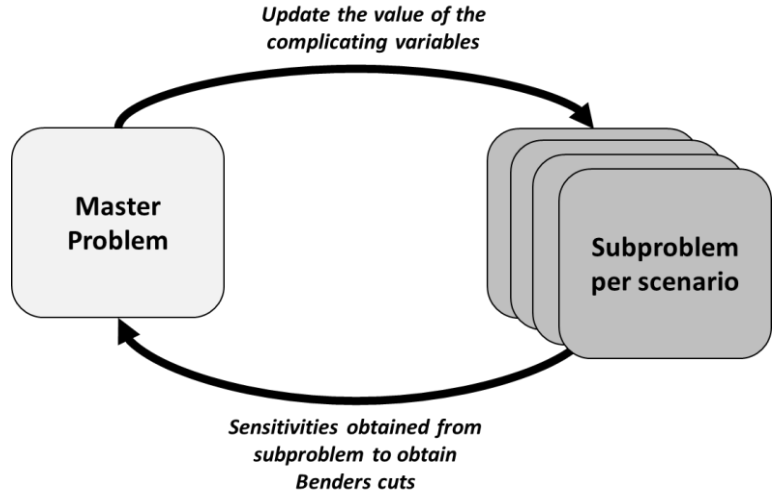

Figure 1: Benders decomposition algorithm.

In the scenario-based decomposition approach, some of the first-stage variables of the two-stage stochastic model are considered as complicating variables. There is only one master problem and a subproblem for each scenario. The subproblem for each scenario $s$ at iteration $\iota$ is given by equations (11)-(20).

$$\underset{\substack{P^{prod,(\iota)}_{l_i,f}, P^{trans,plan,(\iota)}_{l_i,l_j}, P^{trans,actual,(\iota)}_{l_i,l_j} \\ P^{trans,dev(\iota)}_{l_i,l_j,s}, P^{short,(\iota)}_{l_i,s}, P^{excess,(\iota)}_{l_i,s}}}{\text{Minimize}} \sum_{s \in S} \sum_{l_i \in L} \sum_{l_j \in L} \sum_{f \in F} \rho_s \left( C^{prod}_{l_i,f} P^{prod,(\iota)}_{l_i,f,s} \right. \quad (11)$$

$$\left. + \gamma^{short}_{l_i} P^{short,(\iota)}_{l_i,f,s} + \gamma^{trans}_{l_i,l_j} P^{trans,dev,(\iota)}_{l_i,l_j,s} \right)$$

Subject to

$$P^{prod(\iota)}_{l_i,f,s} \leq \bar{P}_{l_i,f} \quad \forall l_i \in L, \forall f \in F \quad (12)$$

$$P_{l_i,f,s}^{prod(\iota)} = \hat{P}_{l_i,f} \quad \forall l_i \in L, \forall f \in F_C \tag{13}$$

$$P_{l_i,f,s}^{prod(\iota)} \leq \hat{P}_{l_i,f} \quad \forall l_i \in L, \forall f \in F_U \tag{14}$$

$$P_{l_i,f,s}^{prod(\iota)} \leq \hat{P}_{l_i,f,s} \quad \forall l_i \in L, \forall f \in F_V \tag{15}$$

$$P_{l_i,l_j,s}^{trans,actual(\iota)} \leq P_{l_j,l_i}^{trans,plan(\iota)} \quad \forall l_i \in L, \forall l_j \in L \tag{16}$$

$$P_{l_i,l_j,s}^{trans,dev(\iota)} = P_{l_j,l_i}^{trans,plan(\iota)} - P_{l_i,l_j,s}^{trans,actual(\iota)} \quad \forall l_i \in L, \forall l_j \in L \tag{17}$$

$$\sum_{f \in F} P_{l_i,f,s}^{prod(\iota)} - \sum_{l_j \in L} P_{l_i,l_j,s}^{trans,actual(\iota)} + \sum_{l_j \in L} P_{l_j,l_i,s}^{trans,actual(\iota)} + P_{l_i,s}^{short(\iota)} - P_{l_i,s}^{excess(\iota)}$$
$$= D_{l_i,s} \quad \forall l_i \in L \tag{18}$$

$$P_{l_i,f,s}^{prod(\iota)}, P_{l_i,l_j,s}^{trans,actual(\iota)}, P_{l_i,s}^{short(\iota)}, P_{l_i,s}^{excess(\iota)} \geq 0 \quad \forall l_i \in L, \forall f \in F \tag{19}$$

$$P_{l_i,l_j}^{trans,plan(\iota)} = P_{l_i,l_j}^{trans,plan,fixed(\iota)} \quad : \mu_s^{trans,plan(\iota)} \quad \forall l_i, l_j \in L \tag{20}$$

The objective function of the subproblem (11) aims to minimize the total scenario-dependent cost. Constraints (12)-(19) encompass all the stochastic constraints relevant to scenario-dependent decisions. In this framework, the complicating variable is $P_{l_i,l_j}^{trans,plan}$, which represents the amount of energy allocated for interchange purposes. Consequently, at each iteration, the value of $P_{l_i,l_j}^{trans,plan(\iota)}$ is fixed for all scenarios as stipulated in constraint (20).

The master problem is given by equations (21)-(25).

$$\underset{P_{l_i,l_j}^{trans,plan,(\iota)}}{\text{Minimize}} \sum_{l_i \in L} \sum_{l_j \in L} C_{l_i,l_j}^{trans} P_{l_i,l_j}^{trans,plan(\iota)} + \alpha^{(\iota)} \tag{21}$$

Subject to

$$\alpha^{(\iota)} \geq \sum_{l_i \in L} \sum_{l_j \in L} \sum_{f \in F} \sum_{s \in S} \rho_s \left( C_{l_i,f}^{prod} P_{l_i,f,s}^{prod(k)} + \gamma_{l_i}^{short} P_{l_i,s}^{short(k)} \right.$$
$$\left. + \gamma_{l_i,l_j}^{trans} P_{l_i,l_j,s}^{trans,dev(k)} \right)$$
$$+ \sum_{s \in S} \sum_{l_i \in L} \sum_{l_j \in L} \mu_s^{trans,plan(k)} \left( P_{l_i,l_j}^{trans,plan(\iota)} - P_{l_i,l_j}^{trans,plan(k)} \right) \tag{22}$$
$$\forall k = 1, \dots, \iota - 1$$

$$\alpha^{(\iota)} \geq \alpha^{down} \tag{23}$$

$$P_{l_i,l_j}^{trans,plan(\iota)} \leq \bar{T}_{l_i,l_j} \quad \forall l_i, l_j \in L \tag{24}$$

$$P_{l_i,l_j}^{trans,plan(\iota)} \geq 0 \quad \forall l_i, l_j \in L \tag{25}$$

The master problem objective (21) is to refine the value of the complicating variable, utilizing the sensitivities determined from the subproblems [41]. Constraint (22) represents Benders optimality cuts, adding a new cut in each iteration. Constraint (23) establishes a lower bound for the solution in the first iteration, which is applicable when no Benders cuts have been added yet. Constraints (24)-(25) are drawn from the stochastic problem, specifically addressing aspects pertaining to the complicating variable.

## 5. A Continental United States (US) Power System Case Study

### 5.1. Selection of Application Area

In order to demonstrate the applicability and robustness of our proposed energy system model, we present an in-depth case study focusing on the United States (US) energy sector. The US energy sector is complex and multifaceted, with a total power generation of approximately 4.1 trillion kWh in 2021, of which renewable sources accounted for 20% [43]. The use of renewable energy is increasing year-on-year in the US [44], with the most prominent source being wind energy. Since 2008, wind energy has formed 64% - 77% of the total year utility from renewable sources in the US [45]. Wind energy is available nationwide, with 90,000 turbines installed in distributed wind applications in all 50 states as of 2022 [46].

Previous studies have focused on optimizing energy systems in specific US regions. For example, Diakov et al. [47] developed a production-cost model for the Western Interchange. Several works have also focused on specific energy sources, such as Brown & Botterud [48], which proposed a zero-carbon national electricity system. However, to the extent of our knowledge, no existing models optimize energy planning across the US, a gap our study aims to fill. We selected the continental US as an application area because it encompasses multiple interconnected regions actively engaged in regional energy interchange, offering an extensive scope for our study.

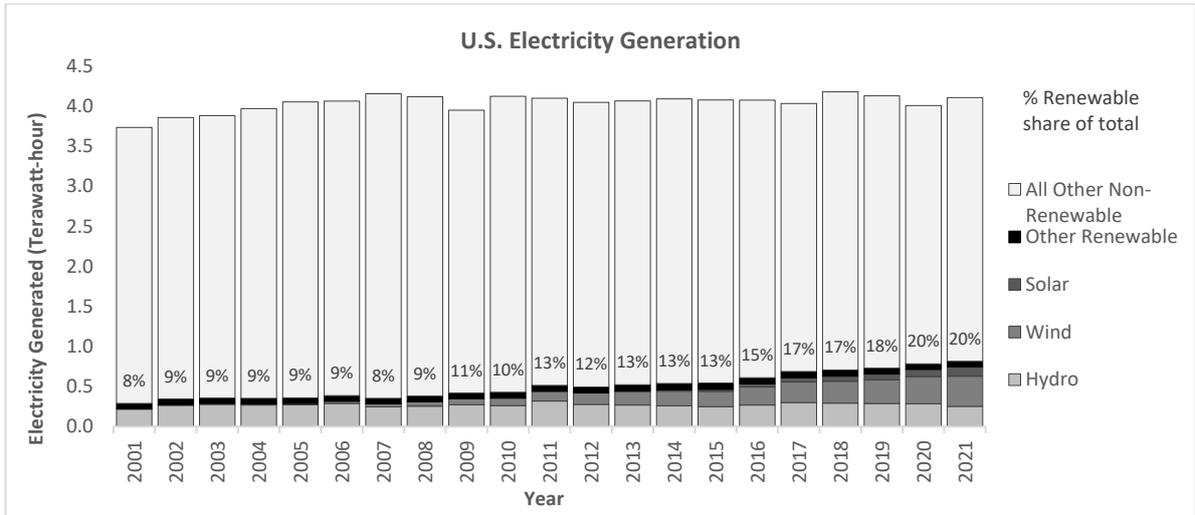

*Figure 2: Amount of electricity generated in the US from 2001-2021.*

### 5.2. Data Sources for the US Application

Data for this work is obtained from the Energy Information Administration (EIA) [49], which is an agency of the US Department of Energy. As shown in Figure 3, the EIA splits up the US into 13 market regions: California (CAL), Carolinas (CAR), Central (CENT), Texas (TEX), Florida (FLA), Mid-Atlantic (MIDA), Midwest (MIDW), New England (NE), New York (NY), Northwest (NW), Southeast (SE), Southwest (SW), and Tennessee (TEN). The types of generators included are wind, solar, hydro, petroleum, natural gas, coal, and nuclear. An "Other" category also exists, which aggregates less common fuel sources, providing a complete overview of energy generation types.

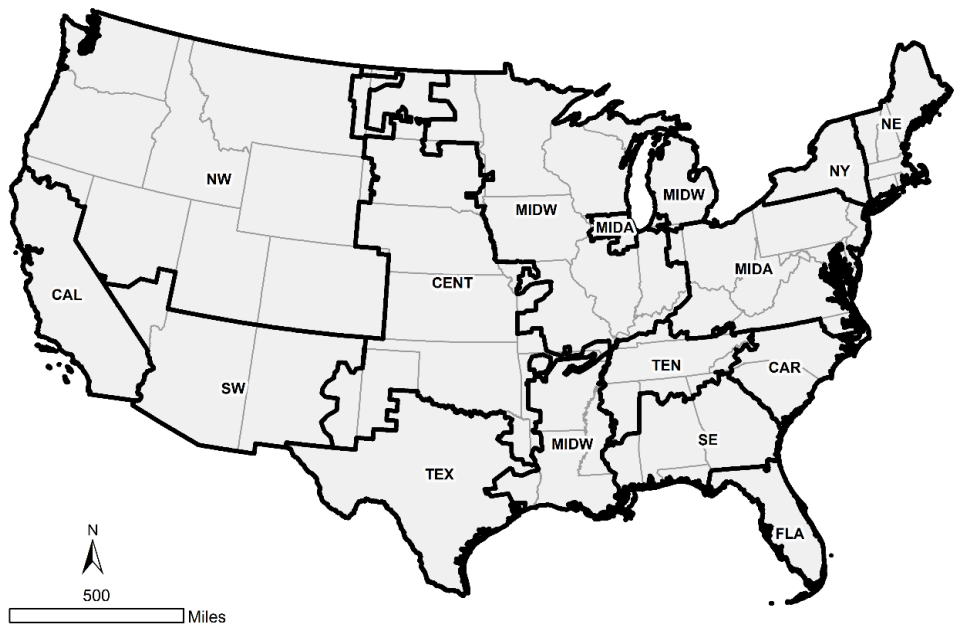

*Figure 3: Map outlining the EIA regions.*

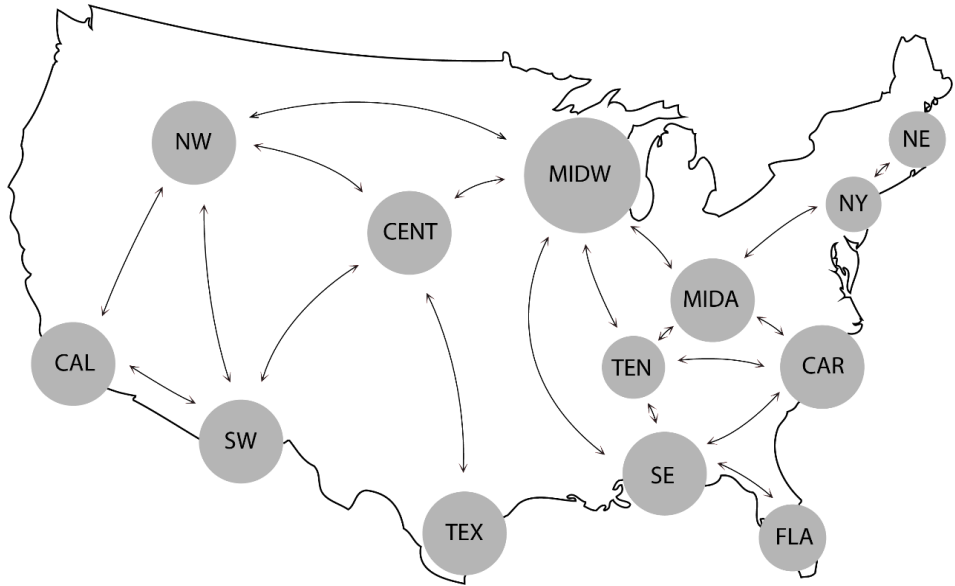

*Figure 4: Map outlining the network of regional energy interchange.*

The following sections describe the scenario-independent and scenario-dependent parameters of the model. Each section includes a description of the parameter, data sources, and notes on any estimates made.

### 5.2.1. Scenario-independent inputs

Parameters that are scenario-independent are those that are assumed to be known to the decision-maker.

#### 5.2.1.1. The cost of generation

In this model, the emphasis is on the cost of electricity generation rather than the price of purchasing electricity. The projected generation costs reported by the International Energy Agency (IEA) [50] are used. Specifically, the median of the levelized cost of electricity (LCOE) for each technology is used.

#### 5.2.1.2. The cost of a shortage

Outage costs are considered reasonably representative measures of system reliability [51], [52]. Some power consumers are extremely sensitive to disruptions, and even a momentary disruption can have devastating effects [53]. Interruption costs can take both direct and indirect forms [54]. There are different methods for estimating outage costs. Woo et al. [55] use the envelope theorem to estimate the non-residential outage cost in the US at $10.37/kWh unserved.

### 5.2.1.3. The cost of transmission

The National Renewable Energy Laboratory (NREL) has created a database of long-distance transmission costs [56]. This work uses the median values for long-distance transmission costs between regions.

### 5.2.1.4. The cost of contractual defaults

Bilateral energy contracts provide price stability and certainty essential for long-term planning [57], a crucial aspect for all parties participating in the agreement. In this work, the electricity transmitted between the different regions is analogized to bilateral contracts, as commonly arranged for large energy consumers. Besides the cost of buying energy, there are potential penalties for deviations from the contracted energy amounts[58]. Due to the lack of specific data on these penalties, a parameter $\kappa^{trans}$ is introduced is the model. This parameter sets the penalty for deviation, $\gamma_{l_i,l_j}^{trans}$, as a product of $\kappa^{trans}$ and the cost of transmission, $C_{l_i,l_j}^{trans}$.

### 5.2.1.5. The capacity of each generating source

The capacity of each generating source was estimated using the EIA dataset [49] by determining the maximum energy produced in the last year for each fuel type in each region for one hour.

### 5.2.1.6. The capacity of transmission

The maximum amount of energy transferred between two locations was estimated using the EIA dataset [49] by calculating the maximum hourly energy transferred between them in the previous year.

### 5.2.2. Scenario-dependent inputs

In stochastic programming, scenarios are developed to account for the uncertainty in the model parameters. Many techniques have been developed for generating scenarios, and each has its strengths and limitations, but no single scenario-generation method is the best for all models [59]. The quality of the scenario generation process affects the accuracy and reliability of the results.

The US Energy Information Administration's (EIA) Hourly Electric Grid Monitor provides up-to-the-hour information on electricity demand across the US electric grid [49]. For parameters that are subject to uncertainty, influential scenarios are developed using a k-means clustering algorithm [60]. This method partitions historical data into $k$ disjoint subsets, or clusters, each represented by a centroid that is used to construct the scenario set. The selection of the optimal number of clusters, $k$, is critical to ensuring the scenarios are representative and capture the variability of the data without overcomplicating the model. To determine the value for $k$, we employ the elbow point method, which involves plotting the within-cluster sum of squares against the number of clusters and identifying the "elbow" point where the rate of decrease sharply changes. This point signifies a balance between model complexity and explanatory power. Once determined, the probability of each chosen scenario is assigned based on the number of data points forming the corresponding cluster.

### 5.2.2.1. Demand Data

Demand data was also obtained from the EIA dataset [49]. All available data at the time was obtained and separated into months by location over all the years. Then, four scenarios for each month in each region were generated using the k-means algorithm. The data was divided into monthly data to include seasonality in the model.

### 5.2.2.2. Available Energy Data

This model separates generators into three categories, as described in Section 3. The generator capacity is used for the first two sets of generators, those with a set amount of energy and those dependent on a controllable resource. For VRRGs, the amount of energy available is scenario-dependent since it relies on environmental factors that are not known ahead of time. Net generation data from the EIA dataset [49] generates scenarios using the k-means algorithm by separating the data into months to include seasonality in the model.

## 5.3. Case Study Results

The presented framework was implemented with MATLAB and modeled with CVX using the Mosek solver. The results of the numerical experiments are presented in the following section, along with their respective implications.

### 5.3.1. Benders convergence

Our numerical experiments found that the Benders decomposition approach is very efficient in finding the optimal solution. On average, it took 4.7 iterations to reach convergence; however, for 78% of the runs, the optimal solution was found within two iterations. The maximum number of iterations to reach convergence among all the experimental runs is 62.

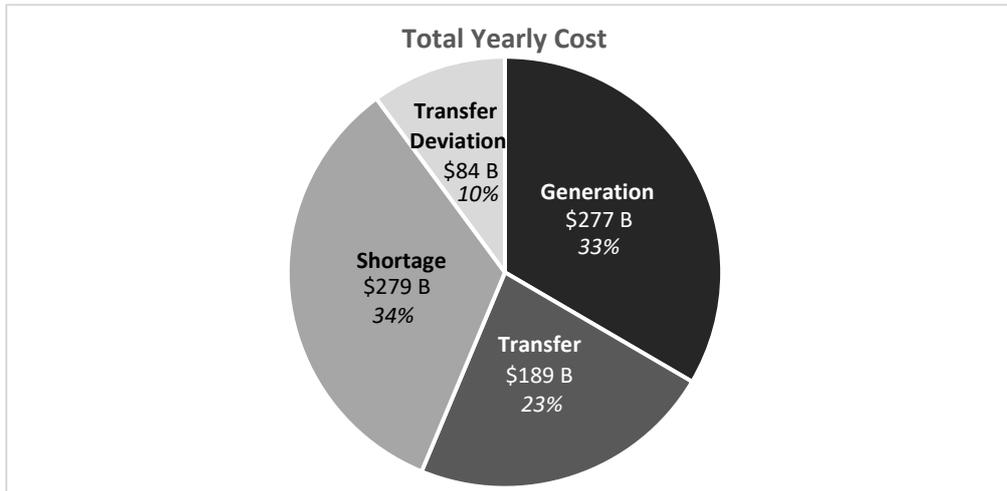

Figure 5: Total yearly cost over all regions split by cost type.

### 5.3.2. Breakdown of costs within the objective function

Section 5.3.2 presents a detailed breakdown of the various cost types within the objective function. The aim is to minimize the total cost, which encompasses generation cost, transfer cost, shortage cost, and penalty for transfer deviation. The framework calculates the total annual cost across all the regions at $829 B (billion). Figure 5 illustrates the proportional contribution of each cost category to this annual total. Table 1 provides key daily value statistics for each cost component. For clarity and consistency, we use the same color coding for different cost categories in the subsequent figures.

|  | Generation Cost ($) | Transfer Cost ($) | Shortage Cost ($) | Transfer Deviation Penalty ($) | Total Cost ($) |
|---|---|---|---|---|---|
| **Minimum** | 26,026,856 | 0 | 0 | 0 | 26,368,267 |
| **Maximum** | 38,703,658 | 117,931,960 | 342,172,025 | 136,343,772 | 634,928,280 |
| **Mean** | 31,664,298 | 21,452,700 | 31,602,206 | 9,414,290 | 94,133,494 |
| **Median** | 31,145,665 | 12,476,517 | 5,695,743 | 0 | 49,309,081 |

Table 1: Key statistics for daily objective function costs.

### 5.3.3. Impact of seasonality

Figure 6 and Figure 7 illustrate the daily results disaggregated by months, revealing the influence of seasonality. Figure 6 shows markedly higher costs in the summer months. This increase aligns with the surge in demand depicted in Figure 7, coupled with the shortfall in energy generation to meet this demand. Seasonal fluctuations in power demand are primarily temperature-driven, with residential and commercial sectors being the most affected [61]. The observed cost rise is attributed to heightened shortage costs, transfer costs, and penalties for transfer deviation, while the generation cost remains constant due to generation capacity constraints. Section 5.3.5 will explore whether augmenting generation capacity could mitigate these costs. Alternative strategies, such as utilization of ESSs, are also considered.

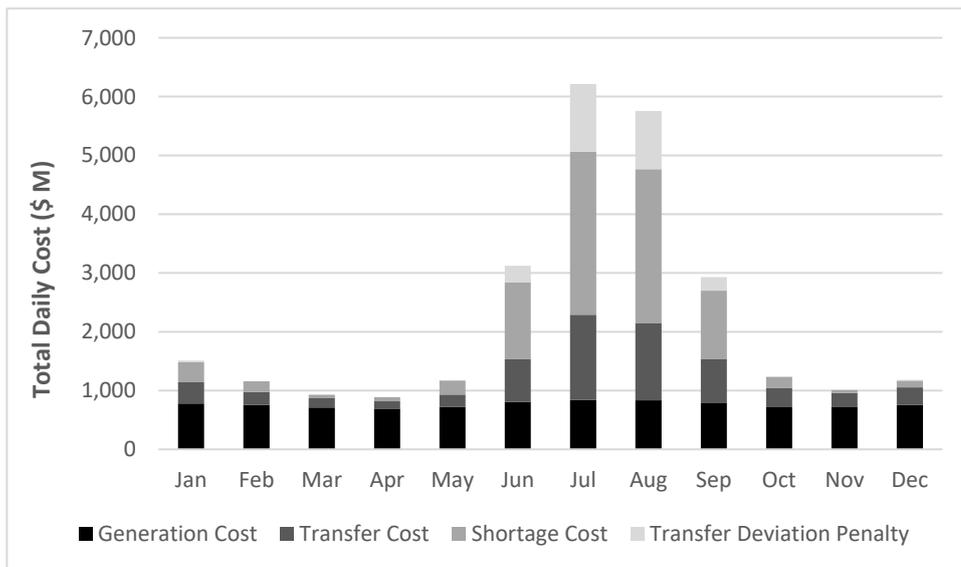

*Figure 6: Comparison of total daily cost for each month over the year.*

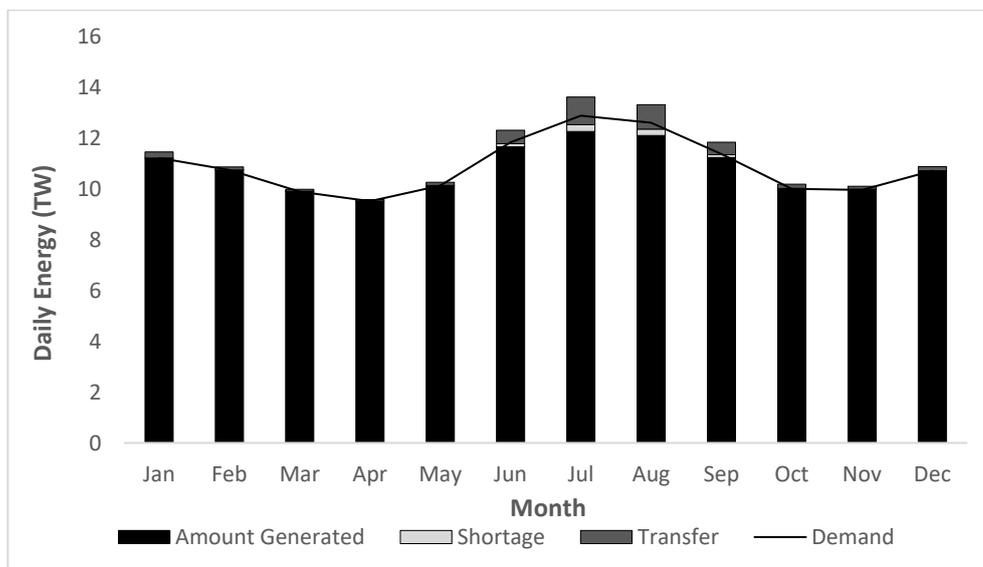

*Figure 7: Total daily energy allocation over all regions by month.*

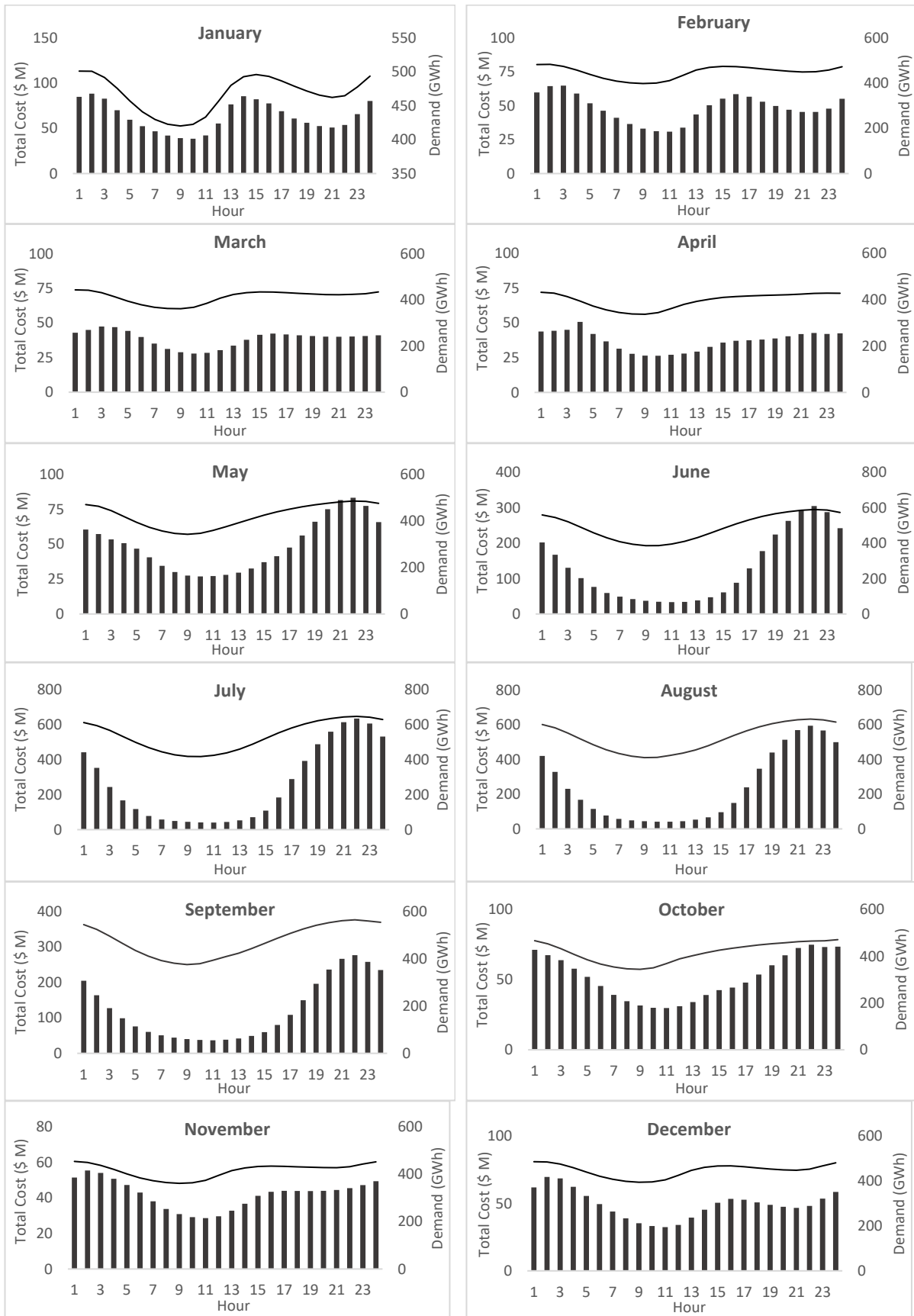

*Figure 8: Total hourly cost for each month in UTC.*

Consideration of national hour-by-hour demand is notably scarce in existing energy generation literature, underlining a significant gap in our understanding of energy patterns. This study contributes substantially to this understudied area by providing detailed insights into the hourly demand variations across the US throughout the year. Figure 8 presents a series of charts illustrating each month's typical daily variation in nationwide cost, displayed in Coordinated Universal Time (UTC). The charts reveal distinct patterns: from May through to October, there is a peak between 5 pm and 5 am, with the highest point around 10 pm. Conversely, two peaks are observed from November through to March –the first around 2 am and the second around 4 pm. Throughout the year, the lowest cost is around 10-11 am. These trends in cost closely align with demand patterns, offering valuable insights for scheduling preventative maintenance, developing demand response strategies, and guiding ESS utilization. In this context, this study serves as a valuable addition to the field, complementing initiatives like the DOE's Demand Response and Storage Integration Study for the Western Interconnection [62].

### 5.3.4. Regional results

Section 5.3.4 delves into the regional analysis of the results. Figure 9 and Figure 10 focus on the location-specific outcomes. Figure 9 reveals that the regions of California (CAL), Mid-Atlantic (MIDA), and Midwest (MIDW) incur the highest total costs. In MIDA and MIDW, a significant proportion of cost arises from generation, correlating with higher energy demand and generation observed for these regions in Figure 10. On the other hand, in the CAL region, the energy generated is insufficient to meet the demand, and the transferred energy is inadequate to mitigate this shortage, resulting in significantly high costs. This finding is corroborated by the 2021 state of emergency declaration in California, which addressed the region's immediate energy deficits [63].

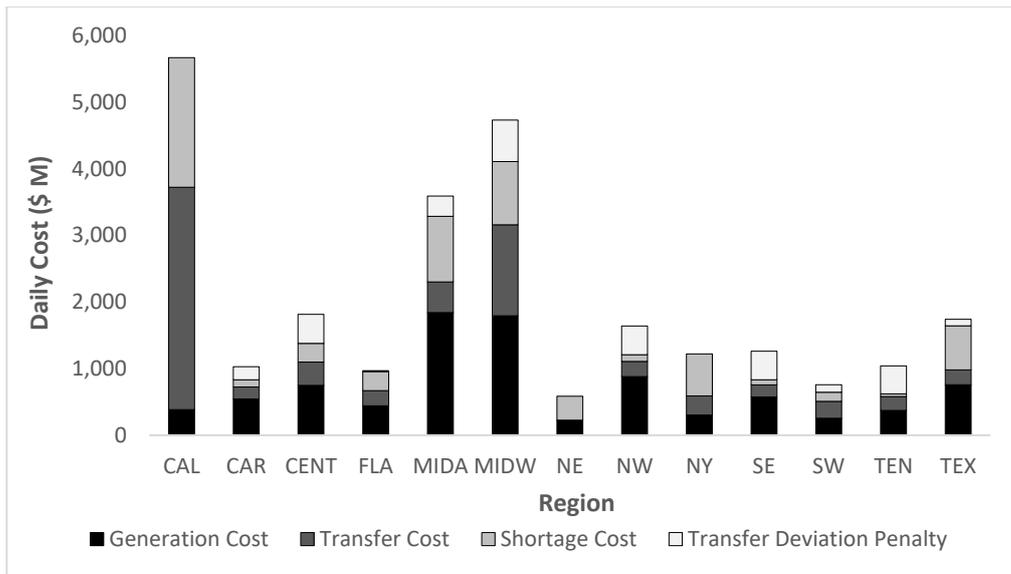

*Figure 9: Breakdown of objective function costs by location*

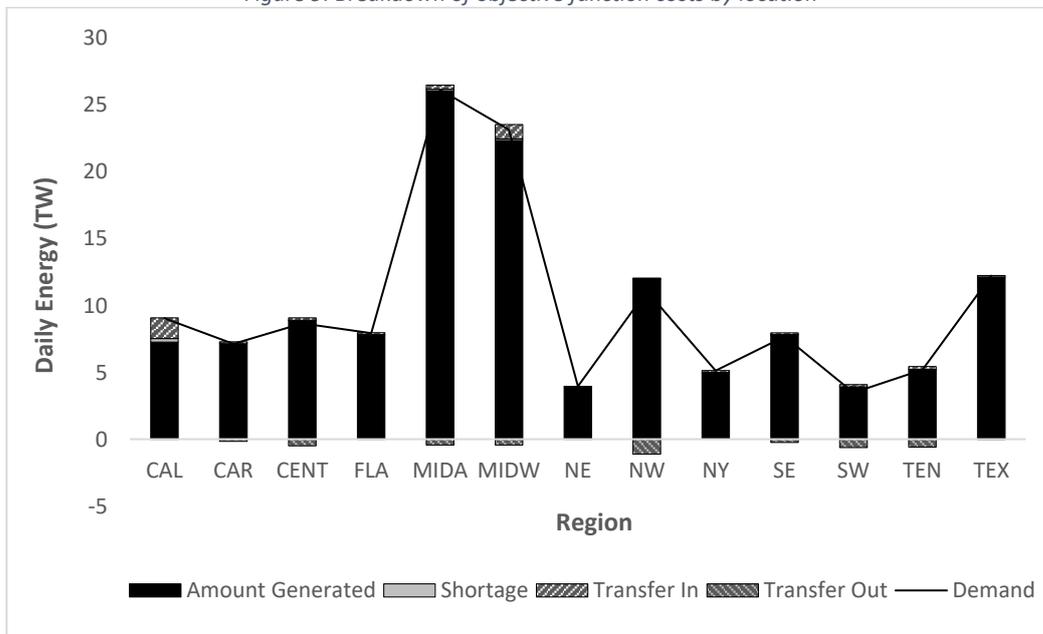

*Figure 10: Average daily energy allocation amount for each region.*

### 5.3.5. Impact of increasing the capacity

Section 5.3.5 examines the potential impact of increased generating capacity on the system's effectiveness. The model was adapted to increase the generating capacity incrementally in a controlled manner to assess how sensitive the total cost is to changes in generating capacity. In each of the 13 iterations, only the generation capacity of a single region was increased by 1 GWh, while the capacities of all other regions were kept at their baseline levels. Figure 11 compares the total nationwide cost between the baseline case and each iteration where one region's capacity is enhanced.

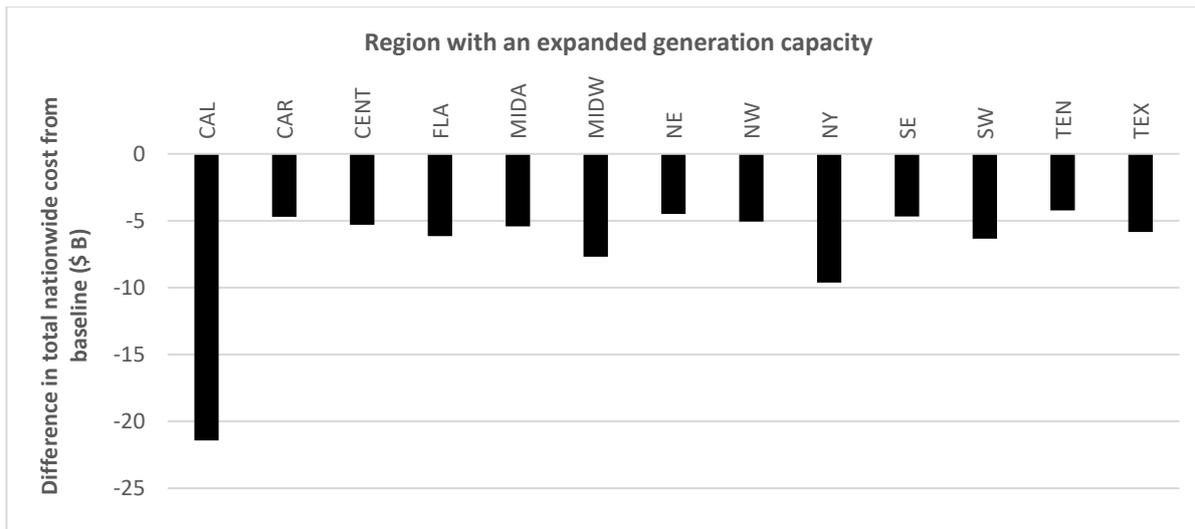

*Figure 11: Result of expanding a region's generation capacity by 1GWh.*

The results show that increasing the generation capacity in any region would decrease the overall nationwide cost. The decrease in cost ranges from $4.2 B to $21.4 B. The most significant overall nationwide cost reduction occurs if the capacity is increased in the CAL region due to the high shortage costs. According to estimates given by the EIA, the utility-scale electric generator construction costs for increasing the capacity of a region by 1GW would be $1.7 B for solar, $1.5 B for wind, and $1.1 B for natural gas [64]. Therefore, increasing the generation capacity in any region would be nationally beneficial.

### 5.3.6. Impact of increasing transmission capacity

Our numerical experiments show that a significant amount of energy shortage can be avoided by sharing energy between locations through regional interchange. However, the potential for energy transmission between regions is bounded by their respective transmission capacities. Numerical experiments were conducted with an assumed infinite transmission capacity to evaluate the impact of the transmission limits. Figure 12 compares the outcomes of this hypothetical scenario with the baseline case.

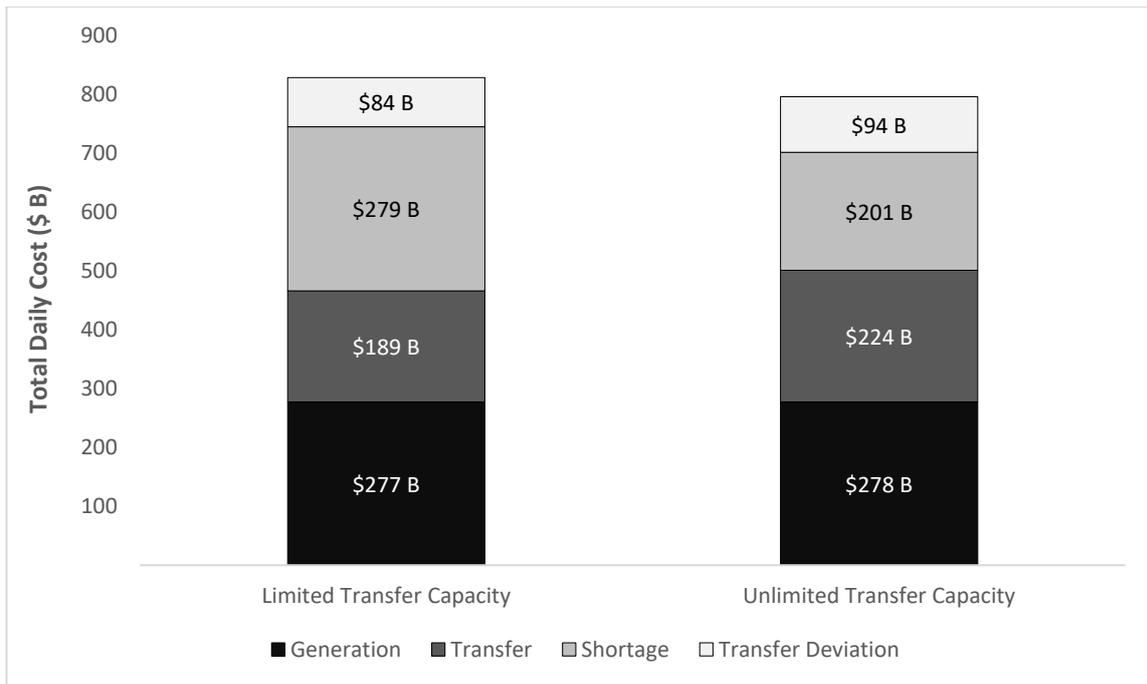

*Figure 12: Comparison of limited and unlimited transfer capacity.*

When comparing the unlimited and limited transfer capacity cases, it is evident that more energy is being generated and transferred in the unlimited transfer capacity case, and the result is lower shortage cost and, thus, a lower total cost. In the unlimited transfer capacity case, the overall cost is decreased by 4%, which is equivalent to $32B yearly. Figure 13 illustrates the average hourly differences between the planned interchange amount between these two cases. Notably, in the unlimited interchange case, more energy is transferred from the NW to the MIDW regions and from the MIDW to the NW, indicating a potential area for transmission line expansion to improve energy distribution efficiency. Additionally, the observed decrease in energy being transferred from TEN to MIDW highlights the redirection of energy due to limited generation capacity.

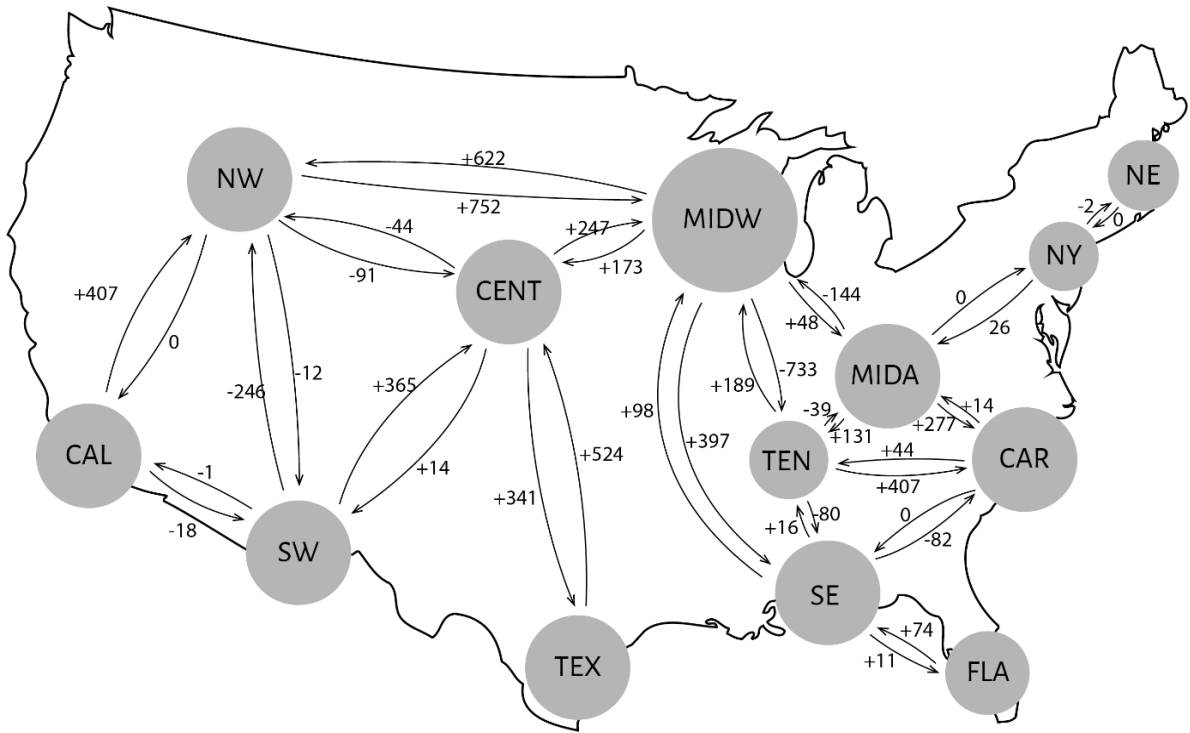

*Figure 13: Average hourly differences in the planned interchange between unlimited and limited generation capacity cases (MWh).*

### 5.3.7. The Expected Value of Perfect Information

The framework presented in this work relies on scenarios to describe the uncertainties at the time of decision-making. The next analysis evaluates the cost of uncertainty to see how much the uncertainty of the scenario-dependent parameters impacts the overall cost and how much the system would benefit from having perfect information. This analysis can also serve as a benchmark to provide a boundary on the gains that can be achieved from the resolution of the uncertainty [65]. One approach to quantifying the cost of uncertainty in decision-making is through a measure in the Value of Information (VOI) framework, known as the Expected Value of Perfect Information (EVPI). EVPI is conceptualized as the maximum price a decision-maker would be willing to pay for completely accurate and certain information. It is calculated by finding the difference between the expected value given perfect information (EV|PI) and the Expected Monetary Value (EMV) under the current state of information.

$$EVPI = EV|PI - EMV \qquad (26)$$

(EV|PI) represents the net benefit of implementing the optimal strategy when perfect information is available, whereas EMV is the average net benefit of the strategies that would be adopted given current, imperfect information [66]. The next two sections present the calculation and comparison of EV|PI and EMV values.

### 5.3.7.1. Expected Value Given Perfect Information

Perfect information is not accessible to generators at the time of decision-making. Each scenario from the stochastic model is run through a deterministic model to approximate the EV|PI, treating the scenario as perfect information. The deterministic model's results represent the realization of a single scenario, whereas the stochastic model provides an expected cost considering all scenarios. Figure 14 plots the results of the deterministic outcomes alongside the results from the stochastic model. Note the large variation in the optimum solution when we use a deterministic modeling approach. This highlights the unsuitability of deterministic modeling for problems where perfect information is unattainable in real-world situations.

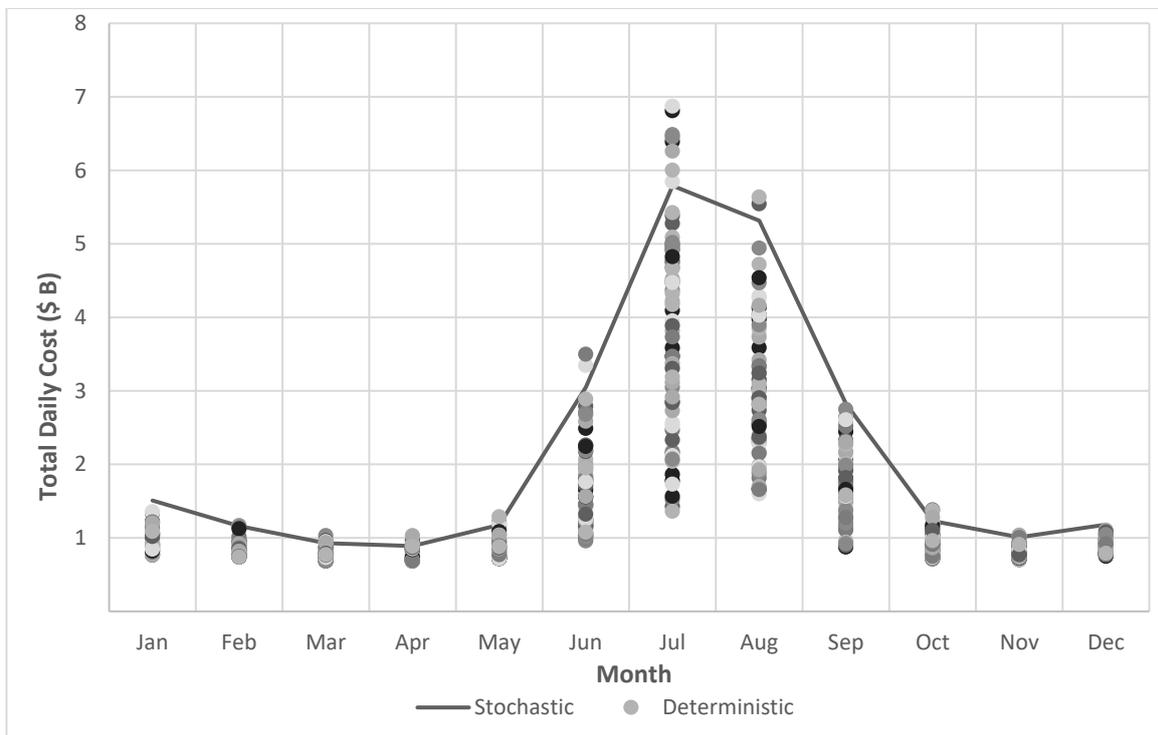

*Figure 14: Comparison of costs using stochastic and deterministic methods ($ MM).*

Deterministic outcomes are used to compute the EV|PI for each month in Table 2, indicating the optimal daily costs under perfect information.

| Expected Value of Total Daily Cost given Perfect Information (EV|PI) ($ M) | | | | | | | | | | | |
|---|---|---|---|---|---|---|---|---|---|---|---|
| Jan | Feb | Mar | Apr | May | Jun | Jul | Aug | Sep | Oct | Nov | Dec |
| 955 | 824 | 757 | 762 | 857 | 1688 | 3877 | 3058 | 1634 | 890 | 793 | 871 |

*Table 2: Total Daily Cost EV|PI.*

### 5.3.7.2. Expected Monetary Value

Given current information, the EMV measures the benefit of an adopted strategy. It is calculated using the expected values of all uncertain parameters in each scenario as inputs to a deterministic model. Table 3 calculates the daily EMV for each month.

| Expected Monetary Value of Total Daily Cost (EMV) ($ M) | | | | | | | | | | | |
|---|---|---|---|---|---|---|---|---|---|---|---|
| Jan | Feb | Mar | Apr | May | Jun | Jul | Aug | Sep | Oct | Nov | Dec |
| 975 | 923 | 812 | 773 | 850 | 1527 | 5107 | 5314 | 1378 | 988 | 713 | 859 |

*Table 3: Total Daily Cost EMV.*

The EVPI is then determined by the difference between the EV|PI and EMV. The result is an average total daily cost of $371 M, translating to an uncertainty cost of 17%. This value represents the theoretical amount a decision-maker should be willing to invest to obtain perfect information, thus quantifying the cost of uncertainty.

## 6. Conclusions and future work

This research work provides a pioneering framework for understanding and optimizing nationwide energy markets and highlights the significance of integrated decision models that can provide optimal resource utilization and cost-efficient strategies for a nation's energy generation systems. Given the variability of renewable sources, such models are crucial for guiding policy decisions, regulatory frameworks, and infrastructure investments. They enable more efficient utilization of available resources, reducing costs and enhancing system reliability.

This study presented a two-stage stochastic programming model to overcome the limitations of deterministic models since perfect information is unavailable when decision-makers make commitments regarding energy generation. Scenarios are generated to handle uncertainty in model parameters, and a scenario-based Benders decomposition approach was adopted to solve the resulting two-stage stochastic optimization model efficiently. This approach provides decision-makers with added flexibility, allowing them to adapt their strategies based on which scenario is realized.

The framework can serve as an invaluable tool to inform policy making, enabling sensitivity analyses and the evaluation of various constraints and policies. A case study focusing on the US energy sector was presented to show the applicability of the proposed framework. Through numerical experiments, significant trends in the total costs of the energy system were uncovered, leading to the identification of implementable policies and strategies to reduce costs for the US energy market.

A key observation of the study is the occurrence of supply-demand gaps in certain regions, with shortage costs accounting for a more significant proportion of the total costs. The shortages mainly occur during summer, aligning with periods of higher energy demand. Analyzing seasonal patterns and integrating national production and demand data enables decision-makers to allocate Energy Storage Systems (ESSs) strategically and schedule preventative maintenance effectively. Furthermore, it was found that increasing the generating capacity in any region would result in an overall cost reduction surpassing the associated construction costs. The current limitations in transmission line capacities impede optimal energy interchange, suggesting that increased investment in both generation and transmission could yield substantial cost savings and better utilization of available resources.

This study is unique in its scope, which is one of the novel features of this research. Still, it is crucial to acknowledge the limitations stemming from the assumptions made in the study, as the existing literature lacks similar integrated models and benchmarks to compare our results against them. Energy planners and policymakers must ensure that the assumptions fit their system and use appropriate datasets before taking any actions or using any policies provided by our framework. Our research shows the necessity and scarcity of accurate region-specific data. Further studies focusing on region-specific data, such as generator capacities and interchange costs, are necessary.

Our study offers a template for future studies and advancements in the critical field of nationwide energy markets. Several promising future research directions exist along which our study can be extended. One avenue involves expanding the scope of the current model to incorporate international entities, such as Mexico's and Canada's energy markets, for the US energy sector. Additionally, the model could be replicated for different collections of regions internationally that currently participate in, or could potentially participate in, energy trade agreements, such as the 'Pan-African' grid [67]. Another potential future research direction would be the inclusion of ESSs in the decision-making. Expanding the model to include additional regions or ESSs would allow the system to participate in market-based arbitrage to satisfy fluctuations in demand. This expansion would enhance the model's applicability and resilience in real-world scenarios.

Regarding the solution approach, the suitability of other decomposition techniques could be investigated to identify opportunities for speeding up the solution process. For instance, multi-cut Benders decomposition can be employed where, in each iteration of the algorithm, a cut is added to the master problem for each scenario instead of adding one cut.

Finally, in the current framework, scenarios are considered independent of one another. Future work could investigate the correlations between different scenarios to develop better forecasts and look at long-term trends for future planning purposes. This would enable more nuanced and informed decision-making for energy planning.

**CRediT authorship contribution statement**

**Razan A. H. Al-Lawati:** Conceptualization, Formal Analysis, Investigation, Mathematical development, Methodology, Validation, Visualization, Writing – original draft. **Tasnim Ibn Faiz:** Conceptualization, Mathematical development, Validation, Writing – review & editing. **Md. Noor-E-Alam:** Conceptualization, Validation, Project administration, Supervision, Writing – review & editing.

**Declaration of Competing Interest**

The authors declare that they have no known competing financial interests or personal relationships that could have appeared to influence the work reported in this paper.

**Acknowledgments**

The authors would like to thank Victoria Okoria Tonbara, an industrial engineering graduate student at the Department of Mechanical and Industrial Engineering, Northeastern University, for her help in the preliminary literature review of this work and proofreading the part of the manuscript.